\documentclass[leqno,10pt]{article}
\usepackage{amsmath}
\usepackage{amssymb,amscd}
\usepackage{latexsym}

\newtheorem{thm}{Theorem}[section]
\newtheorem{lem}[thm]{Lemma}
\newtheorem{prop}[thm]{Proposition}
\newtheorem{cor}[thm]{Corollary}
\newtheorem{rem}[thm]{Remark}
\newtheorem{ex}[thm]{Example}

\newcommand{\RE}{R \!\Join\! E}
\newcommand{\RI}{R\!\Join\! I}
\newcommand{\RM}{R\!\Join\! M}
\newcommand{\VP}{V\!\Join\! P}

\newcommand{\Ker} {\mbox{\rm Ker}}
\newcommand{\Spec} {\mbox{\rm Spec}}

\bigskip

\title{An amalgamated duplication of a ring  \\
 along an ideal: the basic properties 
}
\author{Marco D'Anna\footnote{Partially supported by MIUR, under Grant PRIN 2005-011955.}
\and Marco Fontana\footnote{Partially supported by MIUR, under Grant PRIN 2005-015278.}}

\begin{document}

\maketitle

\hfill{\footnotesize \sl Dedicated to Luigi Salce, on his 60th birthday}%

\bigskip

\begin{abstract}  We introduce  a new general construction,   denoted by $R\!\Join\!E$,
called the amalgamated duplication of a ring $R$ \ along an $R$--module
$E$,  that we assume to be an ideal in some overring of $R$. \ (Note that, when $E^2 =0$,
$R\!\Join\!E$ coincides with the Nagata's idealization
$R\!\ltimes\! E$.)

After discussing the main properties of the amalgamated
duplication $R\!\Join\!E$ in relation with pullback--type
constructions, we restrict our  inve\-sti\-gation to the study of $R\!\Join\!E$ when $E$ is an ideal of $R$.
 Special attention is devoted to the
ideal-theoretic properties of $R\!\Join\!E$ and to the topological
structure of its prime spectrum.
\end{abstract}

\bigskip

%%%%%%%%%%%%%%%%%%%%%%%%%%%%%%%%%%%%%%%%%%%%%%%%%%%%%%%%%%%%%%%%%%%%%%%%%%%%

% SECTION 1

%%%%%%%%%%%%%%%%%%%
\section{Introduction}

If $R$ is a commutative ring with unity and $E$ is an $R$-module,
\sl the idealization $R \!\ltimes\! E$, \rm   introduced by Nagata in 1956 (cf. Nagata's book
\cite{n}, page 2), is a new ring,  containing $R$ as a subring,   where the module $E$ can be
viewed as an ideal such that its square is $(0)$.

This construction has been extensively studied and has many
applications in different contexts (cf. e.g. \cite{r}, \cite{fgr},
\cite{gl}, \cite {h}).
%In particular it has been generalized by Fossum.
% in \cite{fs} with the notion of commutative extension of a
% ring $R$ by an $R$-module.
  Particularly important is the generalization given by Fossum, in \cite{fs}, where he defined \sl  a commutative extension of a ring
$R$ by an $R$--module $E$  \rm to be an exact sequence of abelian groups:
$$
 0 \rightarrow E  \xrightarrow{\iota} S \xrightarrow{\pi} R
\rightarrow 0
 $$
where $S$ is a commutative ring, the map $\pi$ is a ring
homomorphism and the $R$--module structure on $E$ is related to
$S$ and to the maps $\iota$ and $\pi$ by the equation $
 s\cdot \iota(e) = \iota(\pi(s)\cdot e)$ (for all $s\in
S$ and $e \in E$). It is easy to see that the idealization
$R\!\ltimes\! E$ is a very particular commutative extension of $R$
by the $R$--module $E$  (called \sl trivial
extension of $R$ by $E$ \rm in \cite{fs}).

In this paper, we will introduce a new general construction,
called the amalgamated duplication of a ring $R$  along an
$R$--module $E$  (that we assume to be  an ideal in some overring
of $R$ and so $E$ is  an $R$-submodule  of the total ring of
fractions $T(R)$ of $R$) and denoted by $R\!\Join\!E$ (see Lemma
\ref{proj}).

When $E^2 =0$, the new construction $R\!\Join\!E$ coincides with
the idealization $R\!\ltimes\! E$. In general, however, $\RE$ it
is not a commutative extension in the sense of Fossum. One main
difference of this construction, with respect to the idealization
(or with respect to any  Fossum's commutative extension) is that the ring $\RE$ can be a reduced ring (and,  in fact,  it is
always reduced if $R$ is a domain).

 Motivations and some applications of the amalgamated duplication $\RE$ are discussed more in detail in two recent papers  \cite{d'a}, \cite{d'a-f}. More precisely,  M. D'Anna \cite{d'a} has studied some properties  of this
construction in case $E=I$ is a proper ideal of $R$, in order to
construct reduced Gorenstein rings    associated to Cohen-Macaulay rings  and he has applied
this  construction   to curve singularities.   M. D'Anna and M. Fontana in \cite{d'a-f}  have considered the case of the amalgamated duplication of a ring,  in a not necessarily Noetherian setting,  along a multiplicative-canonical ideal in the sense of Heinzer-Huckaba-Papick \cite{hhp}.

The present  paper is devoted  to a more systematic investigation of the ge\-ne\-ral construction $R\!\Join\!E$, with a particular consideration to the
ideal-theoretic properties and to the topological  structure of its prime spectrum.
 More precisely, the  paper is divided in two parts:  in Section 2 we
study the main properties of the amalgamated duplication
$R\!\Join\!E$. In particular  we give a presentation of  this ring as a pullback
(cf. Proposition \ref{pullback}) and from this fact (cf. also \cite{F}, \cite{H-G}) we obtain
several connections between the properties of $R$ and the
properties of $\RE$ and some useful information about Spec($\RE$)
(cf.  Remark   \ref{spec}).

In Section 3 we consider the case when $E=I$ is an ideal of $R$;
this situation allows  us   to deepen the results obtained in Section 2;
in particular we give a complete description of Spec($\RI$)  (cf.
Theorems \ref{spec2} and \ref{loc}).

%%%%%%%%%%%%%%%%%%%%%%%%%%%%%%%%%%%%%%%%%%%%%%%%%%%%%%%%%%%%%%%%%%%%%%%%%%%%
%%   SECTION 2
%%%%%%%%%%%%%%%%%%%%%
\section{The general construction}

In this section we will study the construction of the ring $\RE$
in  a  general setting. More precisely,   $R$ will always be a commutative ring
with unity, $T(R) \ (:= $ $ \{\mbox{regular elements}\}^{-1}R)$ its
total ring of fractions and $E$ an $R$-submodule of $T(R)$.
 Moreover,  in order to construct the ring $\RE$, we are interested in those
$R$-submodules of $T(R)$ such that $E \cdot E \subseteq E$.

%LEMMA 2.1 {alpha}
\begin{lem} \label{alpha} Let $E$ be an $R$-submodule of $T(R)$
and let $J$ be an ideal of $R$.
\begin{description}
  \item[(a)]
$E\cdot E \subseteq E$ if and only if there exists a subring $S$ of
$T(R)$
containing $R$ and $E$, such that $E$ is an ideal of $S$.

 \item[(b)]  If $E\cdot E \subseteq E$ then:
$$
R\!+\!E:=\{z=r+e \in T(R)\ |\ r \in R,\ e\in E\}
$$
 is a subring of $(E:E):=\{z \in T(R)\ |\ zE \subseteq E\}\ (
\subseteq T(R))$, containing $R$ as a subring and $E$ as an ideal.
 \item[(c)]   Assume that $E\cdot E \subseteq E$; the canonical ring
 homomorphism
  $
  \varphi : R \ \hookrightarrow \ R\!+\!E \ \rightarrow \ (R\!+\!E)/E\,, \;
  r \mapsto  r+E\,,
  $
  is surjective and $\Ker(\varphi)= E \cap R$.

   \item[(d)]   Assume that $E\cdot E \subseteq E$; the set
   $J\!+\!E:= \{j+e \ | \ j \in J, \ e \in E\}$
is an ideal of $R\!+\!E$ containing $E$ and $(J\!+\!E) \cap
R=\Ker(R\hookrightarrow R\!+\!E \rightarrow (R\!+\!E)/(J\!+\!E))= J\!+\!(E \cap
R)$.
\end{description}
\end{lem}

\noindent {\bf Proof. } (a) It is clear that the implication
``if'' holds. Conversely, set $S:=(E:E)$. The hypothesis that
$E\cdot E \subseteq E$ implies that $E$ is an ideal of $S$ and
that $S$ is a subring of $T(R)$ containing $R$ as a subring.

\smallskip
\noindent (b) It is obvious that $R\!+\!E$ is an $R$-submodule of
$(E:E)$ containing $R$  and $E$. Moreover, let $r, s \in R$ and
$e, f\in E$,  if $z:=r+e$ and $w:=s+f$ $(\in R\!+\!E)$ then
$zw=rs+(rf+se+ef) \in R\!+\!E$ and $zf = rf +ef \in E$.

\smallskip
\noindent (c)  and (d) are straightforward. \hfill $\Box$

\bigskip
  From now on we will always assume that $E \cdot E \subseteq
E$.

In the $R$-module direct sum $R\oplus E$ we can introduce a
multiplicative structure by setting:
$$
(r,e)(s,f)
 := (rs,rf+se+ef)\,, \; \mbox{where $r,s \in R$
and $e,f \in E$ }.
$$
    We denote by $R {\boldsymbol {\dot {\oplus}}} E$ the
direct sum  $R\oplus E$ endowed  also with the multiplication
defined above.

The following properties are easy to check:
%LEMMA 2.2 {R_E}
\begin{lem} \label{R_E}   With the notation
introduced above, we have:
\begin{description}
  \item[(a)] $R {\boldsymbol {\dot {\oplus}}}  E$ is a ring.
    \item[(b)] The map $j: R {\boldsymbol {\dot {\oplus}}} E
\rightarrow R \times (R\!+\!E)$,
    defined by $(r,e)\mapsto (r,r+e)$, is an injective ring
    homomorphism.
    \item[(c)] The map $i: R  \rightarrow R
{\boldsymbol{\dot {\oplus}}} E$,
    defined by $r\mapsto (r,0)$, is an injective ring
homomorphism.   \hfill $\Box$
\end{description}
\end{lem}

%REMARK 2.3
\begin{rem} \label{idealization} \rm  \bf (a) \rm   With the notation
of Lemma \ref{alpha}, note that if $E= S$ is a subring of $T(R)$
containing as a subring $R$, then $R\!+\!S= S$. Also, if $I$ is an
ideal of $R$, then $R\!+\!I =R$.

 \bf (b) \rm In the statement  of Lemma \ref{alpha} (d), note that,
in general, $J\!+\!E$ does not coincide with the extension of $J$ in
$R\!+\!E$: we have
 $J(R\!+\!E)= \{j+ \alpha  \ |\ j \in J, \ \alpha   \in JE \}\subseteq
J\!+\!E$, but the inclusion can be strict (cf.  Lemma 
\ref{ideals}  (a), (d) and (e)).

\bf (c) \rm  For an \sl arbitrary \rm $R$-module $E$, M. Nagata
introduced in 1955 \it  the idealization of $E$ in $R$,  \rm
denoted here by $R\!\ltimes\! E$, which is the $R$--module
$R\oplus E$ endowed with a multiplicative structure defined by:
$$
(r,e)(s,f) := (rs,rf+se)\,, \; \mbox{ where $r,s \in R$
and $e,f \in E$ }
$$
\noindent (cf. \cite{n-1955} and also Nagata's book
\cite[page 2]{n} and Huckaba's book \cite[Chapter VI, Section 25]{h}).
The idealization $R\!\ltimes\! E$, called also \it the trivial
extension of $R$ by $E$ \rm  \cite{fs},   is a ring such that  the
canonical embedding $R \hookrightarrow R\!\ltimes\! E$, $r \mapsto
(r, 0)$,  defines a subring  of $R \! \ltimes \! E$  isomorphic to
$R$ and the embedding
$E \hookrightarrow R\!\ltimes\! E$, $e \mapsto (0, e)$,  defines an
ideal  $E^\ltimes$ in  $R\!\ltimes\! E$ (isomorphic as an $R$-module
to $E$), which is
nilpotent of index $2$ (i.e.  $E^\ltimes\cdot E^\ltimes =0$).
Therefore, even if $R$ is reduced, the
idealization $R\!\ltimes\! E$ is not a reduced ring, except in the
trivial case for $E=(0)$, since $R\!\ltimes\! (0) = R$.
 Moreover, if $p_R: R\!\ltimes\! E \rightarrow R$ is the canonical
projection (defined by $(r, e) \mapsto r$), then
$$
0 \rightarrow E \rightarrow  R\!\ltimes\! E \xrightarrow {p_{_{\!
R}}}R \rightarrow 0$$
is an exact sequence.

 Note that  the idealization $R\!\ltimes\! E$ coincides with  the
ring   $R {\boldsymbol {\dot {\oplus}}} E$ (Lemma \ref{R_E})  if
and only if $E$ is an $R$-submodule of $T(R)$ that is nilpotent of
index $2$ (i.e. $E\cdot E=  (0)$).
\end{rem}

%LEMMA 2.4 {proj}
\begin{lem}\label{proj}
  With the notation of Lemma \ref{R_E}, note that $\delta :=
j\circ i : R \hookrightarrow R \times (R\!+\!E)$ is the diagonal
embedding and set:

$$ \aligned
R^{\vartriangle}& :=  (j \circ i) (R)= \{(r,r) \ | \ r \in R
\} \ \ \ \text{and}\\
R \Join E &:=  j(R {\boldsymbol {\dot {\oplus}}} E)=
\{(r,r+e) \ | \ r \in R,\
e\in E \}\,.
\endaligned
$$
We have:
\begin{description}
   \item[(a)]    The canonical maps $ R \cong
R^{\vartriangle} \subseteq \RE \subseteq
R \times T(R)$ are ring homomorphisms.

  \item[(b)] $\RE$ is a subdirect product of  the rings  $R
$ and $(R\!+\!E)$,   i.e.  if $\pi_i$ ($i=1,2$) are the projections of
$R \times (R\!+\!E)$
  onto  $R$ and $R\!+\!E$, respectively,  and  if   $\boldsymbol{\mathfrak
O_i}:
=\Ker(\pi_i|_{R\Join E})$,
  then $(\RE)/\boldsymbol{\mathfrak O_1} \cong R$,
$(\RE)/\boldsymbol{\mathfrak O_2} \cong R\!+\!E$ and
$\boldsymbol{\mathfrak O_1}\cap \boldsymbol{\mathfrak O_2}=(0)$.
\end{description}
\end{lem}

\noindent {\bf Proof.}  (a) is obvious.  For   (b) recall that $S$
is  a subdirect product of
a family of rings $\{R_i \ | \ i \in I \}$ if there exists a ring
monomorphism $\varphi : S \hookrightarrow
\prod_i R_i$ such that, for each $i \in I$, $\pi_i\circ \varphi: S
\rightarrow R_i$ is a surjection  (where $\pi_i:  \prod_i
R_i\rightarrow R_i$ is the canonical projection) \cite[page
30]{La}.    Note also
that $\boldsymbol{\mathfrak O_1}= \{(0,e) \ | \ e \in E \}$ and
$\boldsymbol{\mathfrak O_2}=\{(\varepsilon,0) \ | \ \varepsilon \in E
\cap R \}$.
The conclusion is straightforward   (cf. also \cite[Proposition
10]{La}).
 \hfill $\Box$

\medskip
  We will call the ring $\RE$, defined in Lemma \ref{proj}, {\it
the amalgamated duplication of a ring along an $R$ module $E$};
the reason for this name will be clear after studying the prime
spectrum of $\RE$ and comparing it with the prime spectrum of $R$
(see Proposition \ref{spec}).
 The following is an easy consequence of the previous lemma.

%COROLLARY 2.5  {Rdomain}
\begin{cor}\label{Rdomain}  With the notation of Lemma \ref{proj},
 the following properties are equivalent:
\begin{description}
  \item[(i)] $R$ is a domain;
  \item[(ii)] $R\!+\!E$ is a domain;
  \item[(iii)] $ \boldsymbol{\mathfrak O_1} $ is a prime ideal of
$\RE$;
  \item[(iv)] $ \boldsymbol{\mathfrak O_2} $ is a prime ideal of
$\RE$;
  \item[(v)] $\RE$ is a reduced ring and $ \boldsymbol{\mathfrak O_1}
$ and
$\boldsymbol{\mathfrak O_2}$ are prime
  ideals of $\RE$.  \hfill$\Box$
\end{description}
\end{cor}

We will see in a moment that $R$ is a domain if and only if $
\boldsymbol{\mathfrak O_1} $ and $ \boldsymbol{\mathfrak O_2} $
are the only minimal prime ideals $\RE$ (cf. Remark \ref{minp}).

% PROPOSITION 2.6 {pullback}
\begin{prop}\label{pullback}
  Let  $v:  R\times (R\!+\!E)  \twoheadrightarrow  R\times
((R\!+\!E)/E)$ and $u: R \hookrightarrow R\times ((R\!+\!E)/E)$ be  the
natural ring homomorphisms defined, respectively, by $v((x,r+e))
:=(x,r+E)$ and $u(r):=(r,r+E)$, for each $x, r\in R$ and $e \in
E$.
 Then  $v^{-1}(u(R))=\RE$.  Therefore, if
 $v'\ (:=\pi_{1}|_{R\Join E}): \RE
\twoheadrightarrow R$ is the canonical map defined by  $(r,r+e)
\mapsto r$ (cf. Lemma \ref{proj}) and $u': \RE \hookrightarrow
R\times  (R\!+\!E)$ is the natural embedding, then the following diagram:
 $$
\begin{CD}
\RE   @> v'>>    R  \\
@V u'VV       @V u  VV  \\
R\times (R\!+\!E)   @>v>> R\times ((R\!+\!E)/E)
\end{CD}
$$
is a pullback.

 \end{prop}

\noindent {\bf Proof.} Since $E$ is an ideal of $R\!+\!E$
  (Lemma \ref{alpha} (b)),
$ \boldsymbol{\mathfrak O_1} =(0) \times E$ is a common ideal of
$v^{-1}( u(R))$ and $R \times (R\!+\!E)$. Moreover, by definition,
if $x, r\in R$ and $e \in E$, then  $(x,r+e) \in v^{-1}( u(R))$ if
and only if $(x,r+E) \in u(R)$, that is $x-r \in E$.  Therefore we
conclude that  $v^{-1}( u(R) )=\RE$.   The second part of the
statement follows easily from the fact that  $v^{-1}( u(R) )=\RE$
and  $(\RE)/\boldsymbol{\mathfrak O_1} \cong R$, with $
\boldsymbol{\mathfrak O_1}  = \Ker(v')$  (Proposition \ref{proj}
(b)).
 \hfill $\Box$

%COROLLARY 2.7 {fingen}
\begin{cor}\label{fingen}  The ring
$R\times (R\!+\!E)$ is a finitely generated ($\RE$)--module. In
particular,  $\RE \subseteq R \times (R\!+\!E)$ is an integral
extension and $\dim(\RE)=\dim(R \times (R\!+\!E))   = \sup \{\dim(R),
\dim(R\!+\!E)\}$.
\end{cor}

\noindent {\bf Proof.} Clearly $u: R \ \hookrightarrow \ R \times
((R\!+\!E)/E)$ is a finite ring homomorphism, since $R \times
((R\!+\!E)/E)$ is generated by $(1,0)$ and $(0,1)$ as $R$--module.
Since $u$ is finite,  also $ u':  \RE \ (=  v^{-1}( u(R) ))
\hookrightarrow  R \times ((R\!+\!E)/E)$ is a finite ring homomorphism
\cite[Corollary 1.5 (4)]{F}.   Last statement follows from
\cite[Theorems 44 and 48]{Ka}  and from the fact that $\Spec(R
\times (R\!+\!E))$ is homeomorphic to the disjoint union of $\Spec(R)$
and $\Spec(R\!+\!E)$ (cf. also Remark \ref{minp}). \hfill $\Box$

% REMARK 2.8 {minp}
\begin{rem}\label{minp} \rm  Recall that every
ideal of the ring $R \times (R\!+\!E)$ is a direct product of ideals
$I \times J$, with $I$ ideal of $R$ and $J$ ideal of $R\!+\!E$. In
particular,  every prime ideal $Q$ of $R\times (R\!+\!E)$ is  either
of the type $I\times (R\!+\!E)$ or $R \times J$, with $I$ prime ideal
of $R$ and $J$ prime ideal of $(R\!+\!E)$.  Therefore,  in the
situation of Lemma \ref{proj},   if $R$ is an integral domain (and
so $R\!+\!E$ also is an integral domain  by Corollary \ref{Rdomain}),
then $(0) \times (R\!+\!E)$ and $R \times (0)$ are necessarily the
only minimal primes of $R\times (R\!+\!E)$.  By the integrality
property (Corollary \ref{fingen} and \cite[Theorem 46]{Ka}), then
$ \boldsymbol{\mathfrak O_1} =((0) \times (R\!+\!E)) \cap (\RE) =(0)
\times E$ and $ \boldsymbol{\mathfrak O_2} =(R \times (0)) \cap
(\RE) = (R\cap E) \times (0)$ are the only minimal primes of
$\RE$.

Conversely, if $ \boldsymbol{\mathfrak O_1} $ and $
\boldsymbol{\mathfrak O_2} $ are the only minimal primes  of
$\RE$, then clearly  $\RE$ is a reduced ring  (Lemma \ref{proj}
(b)) and, by Corollary \ref{Rdomain}, $R$ is an integral domain.

%COROLLARY  2.9{Noeth}
\begin{cor}\label{Noeth}
  The following statements are equivalent:
\begin{description}
  \item[(i)] $R$  and $R\!+\!E$ are Noetherian;
  \item[(ii)] $R \times (R\!+\!E)$ is Noetherian;
  \item[(iii)] $\RE$ is Noetherian.
  %\item[(iv)] $R\!+\!E$ is Noetherian.
\end{description}
\end{cor}

\noindent{\bf Proof.} Clearly (i) and (ii) are equivalent.  The
statements (ii) and (iii) are equivalent by the Eakin-Nagata
Theorem \cite[Theorem 3.7]{Matsumura}, since  $R\times (R\!+\!E)$ is a
finitely generated ($\RE$)--module (Corollary \ref{fingen}).
\hfill $\Box$

%REMARK 2.10
\begin{rem} \bf (a) \rm  In the situation of Proposition
\ref{pullback},  the pullback degenerates in two cases:

(1)  $v':\RE \rightarrow R$ is an isomorphism if and only if $E = 0$;

(2)
 $u':\RE \rightarrow R \times (R\!+\!E)$ is an isomorphism  if and only if
$E$ is an overring of $R$ (i.e., if and only if $E=R\!+\!E$).

   \bf (b) \rm     By the previous remark, we deduce easily that $R$
Noetherian does not imply in general that $R\!+\!E$ is Noetherian and,
conversely,  $R\!+\!E$ Noetherian    does not imply that $R$ is
Noetherian: take, for instance,  $E$ to be an arbitrary  overring
of $R$.  However, if we assume that $R\!+\!E$ is a finitely generated
$R$-module (cf. also the following Corollary \ref{Noeth+}), then
by the Eakin-Nagata Theorem \cite[Theorem 3.7]{Matsumura} $R$ is
Noetherian if and only if $R\!+\!E$ is Noetherian.

   This same situation described above (i.e. when $E$ is an
arbitrary  overring of $R$) shows that, in Corollary \ref{fingen},
we may have that  $\dim(\RE)=\dim(R)$ or that
$\dim(\RE)=\dim(R\!+\!E)$ (with $\dim(R)\neq \dim(R\!+\!E)$).
 \end{rem}

%COROLLARY  2.11{Noeth+}
\begin{cor}\label{Noeth+}
Assume that $E$ is a fractional ideal of $R$ (i.e. there exists a
regular element $d \in R$ such that $dE \subseteq R$); then the
following statements are equivalent:
\begin{description}
  \item[(i)] $R$  is a Noetherian ring;
  \item[(ii)] $R\!+\!E$ is  a Noetherian $R$-module;
  \item[(iii)] $R \times  (R\!+\!E)$ is a Noetherian ring;
  \item[(iv)] $\RE$ is a Noetherian ring.
  %\item[(iv)] $R\!+\!E$ is Noetherian.
\end{description}
\end{cor}
\noindent{\bf Proof.} By Corollary \ref{Noeth}   and by previous
Remark 2.10 (b), it is sufficient to show that, in this case, $R$
is a Noetherian ring if and only if $R\!+\!E$ is a Noetherian
$R$-module. Clearly, if $R$ is Noetherian, then $E$ is a finitely
generated $R$-module and so $R\!+\!E$ is also a finitely generated
$R$-module and thus it is a Noetherian $R$-module. Conversely,
assume that $R\!+\!E$ is a Noetherian $R$-module; since it is
faithful, by \cite[Theorem 3.5]{Matsumura} it follows that $R$ is
a Noetherian ring. \hfill $\Box$

%COROLLARY 2.12 {intclos}
\begin{cor} \label{intclos}  In the situation described above:
\begin{description}
 \item[(a)]   Let $R'$ and $(R\!+\!E)'$ be the integral closures of $R$
and $R\!+\!E$ in $T(R)$.
  Then $\RE$ and $R\times (R\!+\!E)$ have the same
  integral closure in $T(R) \times T(R)$, which is precisely $R'
\times (R\!+\!E)'$.
  Moreover, if $R\!+\!E$ is a finitely generated $R$-module, then   the
integral closure of $R^{\vartriangle}$  in $T(R) \times T(R)$
(Lemma \ref{proj}) also coincides with $R' \times (R\!+\!E)'$.
  \item[(b)]  If $E\cap R$ contains a regular element, then
  $T(\RE)= T(R\times (R\!+\!E))=T(R) \times T(R)$ and, moreover, $\RE$ and
$R\times (R\!+\!E)$ have the same complete integral closure in $T(R)
\times T(R)$.
\end{description}
\end{cor}

\noindent{\bf Proof.} (a) It is clear that $(x,y) \in T(R) \times
T(R)$ is integral over $R\times (R\!+\!E)$ if and only if $(x,y)
\in R'\times (R\!+\!E)'$. Since the extension    $\RE
\hookrightarrow R\times (R\!+\!E) \ (\subseteq T(R) \times T(R))$
is integral (Corollary \ref{fingen}), we have the first statement.
If, in addition, we assume that   $R\!+\!E$ is a finitely
generated $R$-module, then   the ring  extension $R^{\vartriangle}
\hookrightarrow R \times (R\!+\!E)$ (Lemma \ref{proj}) is finite
(so, in particular, integral) and thus we have the second
statement.

\smallskip
\noindent (b)   Since $E$ is an $R$-submodule of $T(R)$, then
clearly $T(R) = T(R\!+\!E)$, hence it is obvious  that $T(R\times
(R\!+\!E)) = T(R) \times T(R)$. If $e$ is a nonzero regular element of
$E \cap R$, then $(e,e)$ is a nonzero regular element belonging to
$(E\cap R)\times E$, which is a common ideal of $\RE$ and $R\times
(R\!+\!E)$.   From this fact it follows that $\RE$ and $R\times (R\!+\!E)$
have the same total quotient ring \cite[page 326]{Gilmer} and so
$T(\RE)= T(R) \times T(R)$. The last statement follows from
\cite[Lemma 26.5]{Gilmer}. $\mbox{  }$ \hfill $\Box$

\medskip

 Note that,  in Corollary \ref{intclos} (b), the
assumption that $E\cap R$ contains a regular element is essential,
since if $E$ is the ideal $(0)$ of an integral domain $R$ with
quotient field $K$,  then $R \!\Join \! (0) \cong R$ and so $T(R
\!\Join \! (0)) \cong K$, but $T(R \times R) = K \times K$.
\end{rem}

%REMARK  2.13  {spec}
\begin{rem} \label{spec}  \rm   
Using  Theorem 1.4 (c) and Corollary
1.5 (1) of \cite{F},  the previous Proposition \ref{pullback} and Corollary \ref{fingen}  can be used to give a scheme-theoretic description of  Spec$(\RE)$ and Spec$(R
\times (R\!+\!E))$. We do not give here many details, since in the following Section 3 we will prove directly  and in a more elementary way the most part of the statements contained in this remark for the case $E=I$ is an ideal of $R$. 

  Recall that if $f: A \rightarrow B$
is a ring homomorphism, $f^a: \Spec(B) \rightarrow \Spec(A)$  denotes,  as usual,
the continuous map canonically associated to $f$,
i.e. $f^{a}(Q) := f^{-1}(Q)$, for each $Q \in \Spec(B)$; if $I$ is
an ideal of $A$ and if ${\cal X}:= \Spec(A)$, $V_{\cal X}(I)$
  denotes  the Zariski-closed set $\{P\in{\cal X} \mid P
\supseteq I\}$ of ${\cal X}$.

In the situation of Lemma \ref{proj}  and with the notation of
Proposition \ref{pullback},  set $X :=\Spec(R)$, $Y:=\Spec(\RE)$
and $Z:= \Spec(R \times (R\!+\!E))$ and  set $\alpha:= (u')^a: Z
             \rightarrow Y$ and $\beta:=  (v')^a: X
             \rightarrow Y$.  Then   the following properties
hold:
\begin{description}
  \item[(a)] The canonical continuous map $\alpha: Z
             \rightarrow Y$ is surjective.
  \item[(b)] The
             restriction of the map $\alpha: Z \rightarrow Y$
             to $Z \setminus V_Z( \boldsymbol{\mathfrak O_1} )$
gives
rise to a topological homeomorphism:            $$
             \alpha |_{Z \setminus V_Z( \boldsymbol{\mathfrak O_1} )}
:Z \setminus V_Z( \boldsymbol{\mathfrak O_1} ) \
\xrightarrow{\cong} \ Y \setminus V_Y( \boldsymbol{\mathfrak O_1}
)\,.
             $$
         Moreover, for each ${Q} \in
              \Spec(R\times (R\!+\!E))$, with ${ Q} \not\supseteq
\boldsymbol{\mathfrak O_1} $, if $\boldsymbol{{\cal Q}}  : =
\alpha({Q})=  Q \cap (\RE)$,
              then the canonical map $(\RE)_{\boldsymbol{{\cal Q}}}
\rightarrow (R \times (R\!+\!E))_Q$
              is a ring isomorphism.
  \item[(c)]  $\beta: X \rightarrow Y$ defines a canonical
  homeomorphism of $X$ with $V_Y( \boldsymbol{\mathfrak O_1})$;
              moreover, for each $\boldsymbol{{\cal Q}} \in
\Spec(\RE)$ with $ \boldsymbol{{\cal Q}} \supseteq
               \boldsymbol{\mathfrak O_1} $, the canonical ring
homomorphism
$(\RE)/\boldsymbol{{\cal Q}}  \rightarrow  R /v'(\boldsymbol{{\cal
Q}})$
              is an isomorphism.
\end{description}

\end{rem}  

\medskip

  We conclude this section by defining some distinguished ideals of
$\RE$ that are naturally associated to a  given  ideal $J$ of $R$
and  by   giving an example of the general construction.

%PROPOSITION 2.14{J_i}
\begin{prop}\label{J_i} In the situation of Proposition
\ref{pullback} and with the notation of Lemma \ref{alpha}, for
each ideal $J$ of $R$ we can consider the following ideals of
$\RE$:
$$   \boldsymbol{{\cal J}_1} :=v'^{-1}(J)\,,\; \;
\boldsymbol{{\cal J}_2} :=u'^{-1}(R\times J(R\!+\!E)) \;\; \mbox{and}
\; \;    \boldsymbol{{\cal J}_0} :=J^e:=J(\RE)\,.$$ Then we have:
\begin{description}
 \item[(a)] $\boldsymbol{{\cal J}_1} = u'^{-1}(J\times
(R\!+\!E))=u'^{-1}(J\times (J\!+\!E))= \{(j,j+e)\ | \ j \in J, \ e \in E
\}$\,.
\item[(b)]  $\boldsymbol{{\cal J}_0} =\{(j, j+\alpha)\ |\ j \in J,
\ \alpha \in JE \}\,.$
  \item[(c)] $    \boldsymbol{{\cal J}} := \boldsymbol{{\cal
J}_1} \cap  \boldsymbol{{\cal J}_2}=u'^{-1}(J \times J(R\!+\!E))$\,.
  \item[(d)] $\boldsymbol{{\cal J}_0} \subseteq    \boldsymbol{{\cal
J}_1} \cap    \boldsymbol{{\cal J}_2}$\,.
\end{description}
\end{prop}

\noindent {\bf Proof.}   (a) and (b) are straightforward.
Statement (c) is obvious, since $J \times J (R\!+\!E) = (J \times
(R\!+\!E)) \cap (R \times J(R\!+\!E))$.  (d) follows from (c) and from the
fact that $J(\RE)\subseteq u'^{-1}(J(R\times
(R\!+\!E)))=u'^{-1}(J\times J(R\!+\!E))$. \hfill $\Box$

%EXAMPLE 2.15
\begin{ex}
 \rm Let $R := k[t^4,t^6,t^7, t^9]$ (where $k$ is a field and $t$
an indeterminate), $S :=  k[t^2,t^3]$ and $E :=  (t^2,
t^3)S=t^2k[t]$. We have that $R\!+\!E=S$ and hence
$$\aligned
\RE=&\{(f(t),g(t)) \ | \ f \in R,\ g \in  S \ \text{and}\  g-f
\in
E \} = \\
=&\{(f(t),g(t)) \ | \ f \in R,\ g \in  S \ \text{ and}\  f(0)=g(0) \}\ .
\endaligned
$$
Since $E$ is a maximal ideal of  $S$, the prime ideals in $R\times
S $ containing $\boldsymbol{\mathfrak O_1}$ are either of the form
$P \times  S $, for some prime ideal $P$ of $R$, or $R \times E$;
hence the primes not containing $\boldsymbol{\mathfrak O_1}$ are
of the form $R \times Q$, with $Q \in $ Spec$(S)$ and $Q \neq E$.

By  Remark \ref{spec} and  Proposition \ref{J_i},  we have that if $P$ is a
prime in $R$, the ideal $\boldsymbol{\mathcal P_1}=(v')^{-1}
(P)=(u')^{-1}(P \times S)=\{(p,p+e)\ | \ p \in P, \ e \in E \}$
is a prime in $\RE$, containing $\boldsymbol{\mathfrak O_1}$,
and $\RE/\boldsymbol{\mathcal P_1}\cong R/P$.
Moreover, with the notation of Proposition \ref{spec}, in this
way we describe completely $V_Y(\boldsymbol{\mathfrak O_1})$.
Notice also that, if we set $M:= (t^4,t^6,t^7,t^9)R$,
then  the maximal ideals $M \times S$ and $R \times E$ of $R \times S$
have the same trace in $\RE$, i.e.
$(R \times E) \cap (\RE)=\{(r,r+e) \ | \ r \in R \cap E, \ e \in
E\} = (M \times S) \cap (\RE)$.

On the other hand, again by   Remark \ref{spec},    we have that
$Y\setminus V_Y(\boldsymbol{\mathfrak O_1})$ is homeomorphic to $Z
\setminus V_Z(\boldsymbol{\mathfrak O_1})$. Hence the prime ideals
of $\RE$ not containing $\boldsymbol{\mathfrak O_1}$ are of the
form $(R\times Q) \cap (\RE)$, for some prime ideal $Q$ of $S$,
with $Q \neq E$.
%{\bf b)} Let $R=k[t^4,t^7, t^9, t^{10}]$ (where $k$ is a field and
%$t$ an indeterminate), $S=k[t^3,t^4,t^5]$ and
%$E=t^3S=kt^3+t^6k[t]$. We have that $R\!+\!E=k[t^3, t^4]$ and hence
%$$
%\RE=\{(f(t),g(t)) \ | \ f \in R,\ g \in R\!+\!E,\  f(0)=g(0)\
%\text{and} f^{(4)}(0)=g^{(4)}(0) \}\ .
%$$
\end{ex}

%%%%%%%%%%%%%%%%%%%%%%%%%%%%%%%%%%%%%%%%%%%%%%%%%%%%%%%%
%%%%%%%%%%%%%%%%%%%%       SECTION 3   %%%%%%%%%%%%%%%%%
%%%%%%%%%%%%%%%%%%%%%%%%%%%%%%%%%%%%%%%%%%%

\section{The prime spectrum of $R \Join I$}

  In this section we study the case when the $R$-module  $E=I$
is an ideal of $R$ (that we will assume to be proper and different
 from  $(0)$, to avoid the trivial cases); in this situation $R+I=R$.
We start with applying to this case some of the results we
obtained in the general situation.

%PROPOSITION 3.1 {pullback-I}
\begin{prop} \label{pullback-I}
Using the notation of Proposition \ref{pullback}, the
following commutative diagram of canonical ring homomorphisms
$$
\begin{CD}
\RI   @>v'>>    R \\
@Vu'VV       @VuVV  \\
R\times R   @>v>> R\times (R/I)
\end{CD}
$$
is a pullback.
The ideal $ \boldsymbol{\mathfrak O_1} =(0) \times I   = \Ker(v') =
\Ker(v)$
is a common ideal of $\RI$ and $R\times R$, the ideal
$\boldsymbol{\mathfrak O_2}  = \Ker(\RI   \xrightarrow{u'} R \times R
\xrightarrow{\pi_2} R)$
coincides  with   $I \times (0)=
(I \times (0)) \cap(\RI)$ and
$(\RI)/ \boldsymbol{\mathfrak O_i} \cong R$, for $i=1,2$.

  In particular, if $R$ is a domain then $\RI$ is reduced and $
\boldsymbol{\mathfrak O_1} $ and $ \boldsymbol{\mathfrak O_2} $
are the only minimal primes of $\RI$.
\end{prop}

\noindent {\bf Proof.} The first part is an easy consequence of
Lemma \ref{proj} (b) and Proposition \ref{pullback}; the last
statement follows from Corollary \ref{Rdomain}.  \hfill $\Box$

%REMARK 3.2
\begin{rem}
  \rm Note that,  when $I \subseteq R$, then $\RI := \{(r,r+i) \
|\ r \in R, \ i \in I \}= \{ (r+ i,r)  \ |\ r \in R, \ i \in I \}$. It
follows that we can exchange the roles of $\boldsymbol{\mathfrak
O_1} $ and $\boldsymbol{\mathfrak O_2}$ (and that $
\boldsymbol{\mathfrak O_2}$ is  also  a common ideal of  $\RI$ and
$R\times R$).
\end{rem}

If we  specialize to the present  situation Corollary
\ref{fingen}, Corollary \ref{Noeth+} and Corollary \ref{intclos},
then we obtain:

%COROLLARY 3.3{noeth+intclos}
\begin{cor} \label{noeth+intclos} Let $R'$ (respectively, $R^\ast$)
be the integral
closure (respectively, the complete integral closure) of $R$ in
$T(R)$, we have:
\begin{description}

  \item[(a)]   $\dim(\RI) = \dim(R)$.
  \item[(b)]   $R$ is Noetherian if and only if $\RI$ is Noetherian.

    \item[(c)] The integral closure of $R^{\vartriangle}$ and of
$\RI$  in $T(R) \times T(R)$ coincide with $R' \times R'$.

     \item[(d)]  If $I$ contains a regular element, then
$T(\RI) =T(R) \times T(R)$ and the complete integral closure of $\RI$
in $T(R) \times T(R)$ coincide with $R^\ast \times R^\ast$, which is
the complete integral closure of $R\times R$ in $T(R) \times T(R)$.
\end{description}
\end{cor}

%After studying the relation between Spec$(R \times R)$ and
%Spec$(\RI)$, under the  continuous map $(u')^a$, associated the
%canonical embedding $ u':  \RI \hookrightarrow R \times R$, the
%next goal is to investigate directly the relation between
%Spec$(\RI)$ and Spec$(R)$, under the canonical map associated to
%the diagonal embedding  $\delta:  R \hookrightarrow  \RI$, ($r \
%\mapsto (r,r)$).   As above,  we will identify $R$ with its
%isomorphic image $R^{\vartriangle}$ in $\RI$ and we will denote
%the contraction to $R$ of an ideal $\boldsymbol{{\cal H}}$ of
%$\RI$ by $\boldsymbol{{\cal H}} \cap R$ (instead of $\delta^{-1}
%(\boldsymbol{{\cal H}})$).

 The
next goal is to investigate directly the relations among 
Spec($R \times R$), Spec$(\RI)$,  and Spec$(R)$, under the canonical maps associated to natural embeddings, i.e. 
the diagonal embedding  $\delta:  R \hookrightarrow  \RI$, ($r \
\mapsto (r,r)$)  and the inclusion  $\RI \hookrightarrow R \times R $.
 With a slight
abuse of notation, we identify $R$ with its  isomorphic image
$R^{\vartriangle}$ in $\RI\ ( \subseteq R \times R$) under the
diagonal embedding
(Lemma \ref{proj}) and we  denote the contraction to $R$
of an ideal $\boldsymbol{{\cal H}}$ of $\RI$ (or, $H$ of $R\times R$)
by $\boldsymbol{{\cal H}}\cap R$ (or,  by $H \cap R$).

We start with an easy lemma.

%LEMMA 3.4 {ideals}
  \begin{lem}\label{ideals}     With the notation of Proposition \ref{J_i},
let $J$ be an ideal of $R$. Then:
\begin{description}
\item[(a)] $   \boldsymbol{{\cal J}_1} \ (:=v'^{-1}(J))
=u'^{-1}(J\times R)=u'^{-1}(J\times (J\!+\!I))= \{(j,j+i)\ | \ j \in
J, \ i \in I\}    =: J \! \Join \! I$\,.   If $J=I$, then $I \!
\Join \! I \ (= I \times I)$ is a common ideal of $\RI$ and $R
\times R$.
  \item[(b)] $   \boldsymbol{{\cal J}_2} \ (:=u'^{-1}(R\times J))
=
     \{(j+i,j)\ | \ j \in J, \ i \in I\}$\,.
  \item[(c)] $  {\boldsymbol{\cal J}} :=    \boldsymbol{{\cal
J}_1} \cap    \boldsymbol{{\cal J}_2} =u'^{-1}(J \times J)=
    \{(j,j+i')\ | \ j \in J, \ i' \in I\cap J\}=
    \{(j_1,j_2)\ | \ j_1,j_2 \in J, \ j_1-j_2 \in I\}$\,.
  \item[(d)] $  \boldsymbol{{\cal J}_0} \ (:=J(\RI)) =
       \{(j,j+i'')\ |\ j \in J,\ i'' \in JI \}$ {\rm (cf.
       \cite[Lemma 8]{d'a})}.
  \item[(e)] $  \boldsymbol{{\cal J}_0} \subseteq
\boldsymbol{{\cal J}_1} \cap    \boldsymbol{{\cal J}_2} $\,.
  \item[(f)] $   \boldsymbol{{\cal J}_1} =   \boldsymbol{{\cal
J}_2}  \ \Leftrightarrow\ I \subseteq J$\,.
  \item[(g)] $I+J=R \ \Rightarrow \   \boldsymbol{{\cal J}_0}  =
\boldsymbol{{\cal J}_1}  \cap    \boldsymbol{{\cal J}_2} $\,.
  \item[(h)] $   \boldsymbol{{\cal J}_1}  \cap R=
\boldsymbol{{\cal J}_2}  \cap R =   \boldsymbol{{\cal J}_0}  \cap
R =   \boldsymbol{{\cal J}}   \cap R=J$\,.
\end{description}
\end{lem}   

\noindent{\bf Proof.} (a) is a particular case of Proposition
\ref{J_i} (a). The second part is straightforward.

\smallskip
\noindent (b) Let $r \in R$ and $j \in J$;  we have that $(r, j)
\in \RI$ if and only if $(r, j )=(s,s+i)$, for some $s \in R$ and
$i \in I$. Therefore $r=s=j -i$ and $(r, j )=(j+i', j)$ for some
$i' \in I$.

\smallskip
\noindent (c) Let $j_1, j_2 \in J$; we have that $(j_1,j_2) \in
\RI$ if and only if $(j_1,j_2)=(s,s+i)$, for $s \in R$ and $i \in
I$. Therefore $j_1=s$, $j_2=j_1+i$ and $j_2-j_1=i \in I$.

\smallskip
\noindent Statements (d) and (e)  are particular cases of
Proposition \ref{J_i} ((b) and (d)).

\smallskip
\noindent (f) follows easily  from (a) and (b), since:
 $$
    \boldsymbol{{\cal J}_1} =   \boldsymbol{{\cal J}_2} \
\Rightarrow \ J\!+\!I=J \  \Rightarrow \ I \subseteq J \ \Rightarrow \
\boldsymbol{{\cal J}_1} =   \boldsymbol{{\cal J}_2} \ .
$$
(g) is a consequence of (c) and (d), since $J\!+\!I=R$ implies that $J
\cap I=JI$.

\smallskip

\noindent (h)  It is obvious that $   \boldsymbol{{\cal J}_1} \cap
R=J=   \boldsymbol{{\cal J}_2} \cap R$ and hence, by (c) and  (e),
we also have $  \boldsymbol{{\cal J}} \cap R = \boldsymbol{{\cal
J}_0} \cap R=J$. \hfill $\Box$

\medskip
With the help of the previous  lemma we pass to describe the prime spectrum of $\RI$.  In the following, the residue field at
the prime ideal $Q$ of a ring $A$ (i.e. the field $A_Q/QA_Q$) will
be denoted by $\boldsymbol{k}_A(Q)$.  Part of the next theorem is
contained in \cite[Proposition 5]{d'a}.

%THEOREM 3.5 {spec2}
\begin{thm}\label{spec2}

\begin{description}
\item[ \bf (1)] \it  Let $P$ be a prime ideal of $R$ and consider the ideals
$\boldsymbol{{\cal P}_1}$, $ \boldsymbol{{\cal P}_2} $,
$\boldsymbol{{\cal P}_0}$ and $ {\boldsymbol{\cal P}}$ of $\RI$
as in  Lemma \ref{ideals}  (with $P=J$). Then:

{\begin{description}
   \item[(1, a)]   $\boldsymbol{{\cal P}_1} $ and $
\boldsymbol{{\cal P}_2} $ are the   only   prime ideals of $\RI$
lying over $P$\,.
   \item[(1, b)]    If $P \supseteq I$, then $       \boldsymbol{{\cal
P}_1}  =    \boldsymbol{{\cal P}_2} ={ \boldsymbol{\cal P}}
=\sqrt{\boldsymbol{{\cal P}_0}}   = P\!\Join \!I$.
   Moreover,
   $ \boldsymbol k_R(P)\cong \boldsymbol k_{R\Join I}(
{\boldsymbol{\cal P}} )$\,.
  \item[(1, c)]  If $P \not\supseteq I$ then $       \boldsymbol{{\cal
P}_1}  \neq     \boldsymbol{{\cal P}_2} $. Moreover $
\boldsymbol{{\cal P}} =\sqrt{\boldsymbol{{\cal P}_0}}$ and
$\boldsymbol k_R(P) \cong \boldsymbol k_{R \Join I}(
\boldsymbol{{\cal P}_1}) \cong \boldsymbol k_{R \Join I}(
\boldsymbol{{\cal P}_2})
       $\,.
   \item[(1, d)]  If $P$ is a maximal ideal of $R$ then $
\boldsymbol{{\cal P}_1}$ and $    \boldsymbol{{\cal P}_2}$ are
maximal
   ideals of $\RI$\,.
   \item[(1, e)]  If $R$ is a local ring with maximal ideal $M$ then $\RI$
is a local
  ring with ma\-xi\-mal ideal ${\boldsymbol{\cal M}} =
  \sqrt{ \boldsymbol{{\cal M}_0}} = M \! \Join \! I$
  (using again the notation of  Lemma
\ref{ideals}   for $M=J$) .
%%%%%%
 \item[(1, f)]  $R$ is reduced if and only if $\RI$ is reduced.
\end{description}}

 \item[\bf (2)]  \it    Let $ \boldsymbol{\boldsymbol{{\cal Q}} }$ be a prime ideal of $\RI$ and let $ \boldsymbol{\mathfrak O_1}$ be as in  Proposition \ref{pullback-I}.  Two cases are possible either
$\boldsymbol{{\cal Q}}
\nsupseteq  \boldsymbol{\mathfrak O_1}$ or
$\boldsymbol{{\cal Q}}  \supseteq  \boldsymbol{\mathfrak O_1} $. 
{\begin{description}

\item[(2, a)] If $\boldsymbol{{\cal Q}} \nsupseteq
\boldsymbol{\mathfrak O_1} $, then 
 there exists a unique prime ideal ${Q}$ of $R\times R$ such that $\boldsymbol{{\cal Q}}={Q} \cap(\RI)$ with
${Q} =R \times P$, where $P := \boldsymbol{{\cal Q}}  \cap R$  (and $P \nsupseteq I$).  In this case, with the notation of the previous part  (1),  $\boldsymbol{{\cal
P}_1}  \neq     \boldsymbol{{\cal P}_2} $ and 
$$
\boldsymbol{{\cal Q}}= \boldsymbol{{\cal P}_2} = \{(p+i,p) \ | \ p \in P, \ i \in I\}\ .
$$
  Furthermore,  the canonical ring homomorphisms $ \RI
\hookrightarrow R \times R \xrightarrow{\pi_2} R $ induce for the
localizations the following isomorphisms:
$$
(\RI)_{\boldsymbol{{\cal Q}} } \cong (R \times R)_{{Q}}= (R \times
R)_{R \times P}
\cong R_{P}  \  \mbox{ (thus $ \boldsymbol{k}_{\!R\Join I}(\boldsymbol{{\cal Q}}) \cong \boldsymbol{k}_{R}(P)$)}\,.
$$

\item[(2, b)]  If $\boldsymbol{{\cal Q}}  \supseteq
\boldsymbol{\mathfrak O_1} $, then  there exists a unique prime
ideal $P$ of $R$ such that $\boldsymbol{{\cal Q}} =v'^{-1}(P)$
(or, equivalently, $P = v'(\boldsymbol{{\cal Q}}))$.  With the notation of the previous part  (1), if $P\supseteq I$ then 
$\boldsymbol{{\cal Q}} =  \boldsymbol{{\cal P}_1}=  \boldsymbol{{\cal P}_2}$.
On the other hand, if $P\nsupseteq I$ then
$\boldsymbol{{\cal Q}}  =  \boldsymbol{{\cal P}_1} \ ( \neq  \boldsymbol{{\cal P}_2})$. In both cases,
$$
\boldsymbol{{\cal Q}} = \{(p,p+i) \ | \ p \in P, \ i \in I\}\,.
$$
Furthermore,  the canonical ring homomorphism $v':  \RI
\rightarrow R $ induces the following isomorphism:
$$(\RI)/\boldsymbol{{\cal Q}}
\cong R/P \  \mbox{ (thus $ \boldsymbol{k}_{\!R\Join I}(\boldsymbol{{\cal Q}}) \cong \boldsymbol{k}_{R}(P)$)}\,.$$

\end{description}} 
\end{description}
\end{thm}

\noindent{\bf Proof.}  Note that the composition of the diagonal
embedding $\delta: R \hookrightarrow \RI$, \ ($r \mapsto (r, r)$),
with the inclusion $\RI \subseteq R\times R$, \ ($(r, r+i) \mapsto
(r, r+i)$), coincides with the diagonal embedding $R
\hookrightarrow R\times R$, \ ($r \mapsto (r, r)$), which is a
finite ring homomorphism. Thus,   in particular,  both $R
\hookrightarrow \RI$ and $\RI \subseteq R \times R$ are integral
homomorphisms. Note also that if $Q$ is a prime ideal of $R \times
R$ lying over $P$, then necessarily $Q \in \{P \times R, R \times
P \}$    (Remark \ref{minp}).

\smallskip
\noindent   (1, a)  Note that $       \boldsymbol{{\cal P}_1} =u'^{-1}(P
\times R)$ and $  \boldsymbol{{\cal P}_2}=u'^{-1}(R \times P)$
 (Lemma \ref{ideals});   hence $ \boldsymbol{{\cal P}_1} $ and
$ \boldsymbol{{\cal P}_2}$ are prime ideals lying over $P$. By
integrality, if ${\boldsymbol{\cal Q}} \in \text{Spec}(\RI)$ and
${\boldsymbol{\cal Q}} \cap R=P$, then there exists $\overline Q
\in \text{Spec}(R \times R)$ such that $\overline Q \cap
(\RI)={\boldsymbol{\cal Q}}$ and thus $\overline Q \cap R=P$.
Therefore $\overline Q \in \{P \times R, R \times P \}$ and so
${\boldsymbol{\cal Q}} \in \{       \boldsymbol{{\cal P}_1},
\boldsymbol{{\cal P}_2}\}$.

\smallskip
\noindent  (1, b)   We know already by  Lemma  \ref{ideals}   (f) and
(c) that, if $P\supseteq I$, then $     \boldsymbol{{\cal P}_1} =
\boldsymbol{{\cal P}_2} =  \boldsymbol{{\cal P}} $, hence by part
 (1, a)  we conclude easily that  $\boldsymbol{{\cal P}} =\sqrt{
\boldsymbol{{\cal P}_0}}$. Moreover we have the following sequence
of canonical homomorphisms:
$$
\frac{R}{P} \subseteq \frac{\RI}{\sqrt{\boldsymbol{{\cal P}_0}}} =
\frac{\RI}{\boldsymbol{{\cal P}} } \subseteq \frac{R \times R }{P
\times R} \cong \frac{R}{P} \cong \frac{R \times R }{R \times P} \
,
$$
  from which we deduce the last part of the statement.

\noindent  (1, c)   By  Lemma \ref{ideals}  (e) and (f) we know
that, if $P\not\supseteq I$, then $      \boldsymbol{{\cal P}_1}
\neq     \boldsymbol{{\cal P}_2}$ and $\boldsymbol{{\cal P}_0}
\subseteq   {\boldsymbol{\cal P}} = \boldsymbol{{\cal P}_1} \cap
\boldsymbol{{\cal P}_2}$. By part  (1, a)  and by the integrality of $R
\hookrightarrow \RI$, we conclude easily that $\boldsymbol{{\cal
P}} =\sqrt{\boldsymbol{{\cal P}_0}}$. Finally, as in part  (1, b),  it
is easy to see that $\boldsymbol k_R(P)\cong \boldsymbol k_{R\Join
I}( \boldsymbol{{\cal P}_1} ) \cong \boldsymbol k_{R\Join I}(
\boldsymbol{{\cal P}_2})$.

\smallskip
\noindent  (1, d)  follows by the integrality of $R \subseteq \RI$.

\smallskip
\noindent   (1, e)  follows immediately by part  (1, d)  and part  (1, b).  

\smallskip
\noindent  (1, f)  follows by integrality of $R \hookrightarrow \RI$
and $\RI \subseteq R \times R$ and from the fact that $R$ is
reduced if and only if $R\times R$ is reduced. 

\noindent   (2)  If $P =  {\boldsymbol{\cal Q}}  \cap R$, then  necessarily ${\boldsymbol{\cal Q}}  \in  \{ \boldsymbol{{\cal P}_1}, 
\boldsymbol{{\cal P}_2}  \}$ by (1, a).

\noindent  (2, a)   Since ${\boldsymbol{\cal Q}} \nsupseteq
\boldsymbol{\mathfrak O_1} $,  then   ${\boldsymbol{\cal Q}}  = 
\boldsymbol{{\cal P}_2}$, because $\boldsymbol{{\cal P}_1}  \supseteq \boldsymbol{\mathfrak O_1} $. 
%and  $\boldsymbol{{\cal P}_2} \supseteq  \boldsymbol{\mathfrak O_2} $
 %hence $\boldsymbol{{\cal P}_2}  \nsupseteq \boldsymbol{\mathfrak O_1} $.
%Hence ${Q}=R \times
%P$ for some prime $P$ of $R$ such that $P \nsupseteq I$).
  Note that $\boldsymbol{{\cal P}_2} = (R \times P) \cap \RI$;  it is easy to see that $Q:= R \times P$ is the unique prime of $R \times R$ contracting over ${\boldsymbol{\cal Q}} $.  The elementwise description of $\boldsymbol{{\cal P}_2} $ is a particular case of Lemma  \ref{ideals}  (b).  Last statement follows from the following canonical  inclusions of localizations $R_P \hookrightarrow (\RI)_{\boldsymbol{{\cal Q}} }\hookrightarrow (R \times R)_{{Q}}= (R \times
R)_{R \times P}
\cong R_{P}\,.$

\noindent (2, b)  The first and the last statements are trivial consequences of the fact that $v'$ induces an isomorphism between $\RI/\boldsymbol{\mathfrak O_1}$ and $R$.  It is easy to see that the prime $P$ is such that $P = \boldsymbol{{\cal Q}} \cap R$. Therefore the second statement follows from (1, b). If $P \nsupseteq I$ (and $\boldsymbol{{\cal Q}}  \supseteq
\boldsymbol{\mathfrak O_1} $) then
$\boldsymbol{{\cal Q}}  =  \boldsymbol{{\cal P}_1} \ ( \neq  \boldsymbol{{\cal P}_2})$, since $\boldsymbol{{\cal Q}}$ does  not contain $\boldsymbol{\mathfrak O_2} $ (note that a prime ideal of $\RI$ containing both  $\boldsymbol{\mathfrak O_1} $ and  $\boldsymbol{\mathfrak O_2} $  has a trace in $R$ containing $I$).  
The elementwise description of $\boldsymbol{{\cal P}_1} $ is a particular case of Lemma  \ref{ideals}  (a).
%  Proposition \ref{spec}, there exists a unique prime ideal ${Q}$
%of $R\times R$ such that $\boldsymbol{{\cal Q}}={Q} \cap(\RI)$ and
%${Q}
%\nsupseteq (0) \times I$. Therefore ${Q} \nsupseteq (0) \times
%R$ and so ${Q} \supseteq R \times (0)$. Hence ${Q}=R \times
%P$ for some prime $P$ of $R$ such that $P \nsupseteq I$. 
\hfill $\Box$

% REMARK 3.6
\begin{rem} \rm
In the situation of  Theorem \ref{spec2}, note that, if $P$ is a
prime ideal of $R$, then
 by integrality of $ R \hookrightarrow \RI \subseteq R \times R$,
inside the ring $R\times R$, the prime ideals $P \times R$ and $R
\times P$ are the only minimal prime ideals of $P \times P=
\boldsymbol{{\cal P}_0}(R \times R)=P(R \times R)$, and so
$$
\boldsymbol{{\cal P}_0}(R \times R) = P \times P= (P \times R)\cap
(R \times P)  = \sqrt{\boldsymbol{{\cal P}_0}(R \times R)}
$$
 is a radical ideal of $R \times R$, with
 $$
 (P \times P) \cap (\RI) =  ((P \times R) \cap   (R \times P)  )\cap
(\RI) =  \boldsymbol{{\cal P}_1} \cap  \boldsymbol{{\cal P}_2} =
\boldsymbol{{\cal P}}\,.
 $$
Next example shows that in $\RI$,  in general, $\boldsymbol{{\cal
P}_0}$ is not a radical ideal (i.e. it may happen that
$\boldsymbol{{\cal P}_0} \subsetneq \sqrt{\boldsymbol{{\cal P}_0}}
 = \boldsymbol{{\cal P}}$).
\end{rem}

%EXAMPLE 3.7
\begin{ex}\label{ex1} \rm
Let $V$ be a valuation domain with a nonzero non maximal non
idempotent prime ideal $P$.  (An explicit example can be
constructed as follows: let $k$ be a field and let $X, Y$ be two
indeterminates over $k$, then take $V:= k[X]_{(X)}+Yk(X)[Y]_{(Y)}$
and $P:=Yk(X)[Y]_{(Y)}$. It is well known that $V$ is discrete
valuation domain of dimension 2, and $P$ is the height 1 prime
ideal of $V$ \cite[(11.4), page 35]{n}, \cite[page 192]{Gilmer}.)

In this situation, it is easy to see that   the ideal   $P \times
P$ is a common (radical) ideal of $\VP$ and of its overring $V
\times V$.  Moreover, note that  $\boldsymbol{{\cal P}_0}=P(\VP) =
\{(p, p+x)\ |\ p \in P,\ x \in P^2\}$  (Lemma \ref{ideals}
(d))  and  that $P(V\times V) = P\times P \subset \VP$. More
precisely, by  Lemma  \ref{ideals}  (c), we have:
       $$
       \begin{array}{rl}
        P\times P  =& (P\times P ) \cap   (\VP) =  (P\times V ) \cap
(V \times P) \cap (\VP) \\
        =&  \boldsymbol{{\cal P}_1} \cap  \boldsymbol{{\cal P}_2} =
\boldsymbol{{\cal P}} =     \{(p,p+y)\ | \ p \in p, \ y\in P\cap P =
P\}\,.
        \end{array}
        $$
Clearly, since $P^2 \neq P$, then  $\boldsymbol{{\cal P}_0}
\subsetneq \boldsymbol{{\cal P}}$; for instance if $z\in
P\setminus P^2$, then $(p, p+z) \in \boldsymbol{{\cal P}}
\setminus P(\VP)$.
\end{ex}

  We complete now the description of the affine scheme
$\Spec(\RI)$, initiated in Theorem \ref{spec2},  determining  in particular   the localizations of  $\RI$ in each of
its prime ideals. Part of the next theorem is contained in
\cite[Proposition 7]{d'a}.

%THEOREM 3.8
\begin{thm}\label{loc}
%In the situation of Proposition ref{pullback-I},
  Let $X:=\Spec(R)$,  $Y:=\Spec (\RI)$ and $Z:=\Spec(R\times
R)\cong \Spec(R) \amalg\Spec(R)$ and let $\alpha: Z
\twoheadrightarrow Y $ and $\gamma:  Y \twoheadrightarrow X$ be
the canonical surjective maps associated to the integral
embeddings $\RI \hookrightarrow R\times R$ and $R\cong
R^{\vartriangle} \hookrightarrow\RI$ (proof of Theorem
\ref{spec2}).
\begin{description}

\item[(a)] The restrictions of $\alpha$
  $$
   \alpha \left|_{Z\setminus V_Z(\boldsymbol{\mathfrak O_i})}:
   Z\setminus V_Z(\boldsymbol{\mathfrak O_i}) \ \longrightarrow \
   Y \setminus V_Y(\boldsymbol{\mathfrak O_i})\right.
  $$
(for $i=1,2$) are scheme isomorphisms, and  clearly
  $$
  Z\setminus V_Z(\boldsymbol{\mathfrak O_i}) \cong
 %\left( (X \amalg X)\setminus (X \amalg V_{X}(I))\right)=
  X \setminus V_{X}(I)\ .
$$
In particular, for each prime ideal $P$ of $R$, such that $P
\not\supseteq I$, if we set $\overline
  {P_1}:=P\times R$ and $\overline{P_2} :=R \times P$ we have
  $  \boldsymbol{\cal {P}}_{i} := \overline P_i \cap (\RI)$, for
$1 \leq i \leq 2$ and
%$\boldsymbol{\cal {P}}_{1} $ and
%$\boldsymbol{\cal {P}}_{2}$ are distinct prime ideal of $\RI$ and
%they are the only prime ideals of $\RI$ contracting onto $P$.
%Moreover,
the following canonical ring homomorphisms are isomorphisms:
  $$
     R_P \longrightarrow (\RI)_{  \boldsymbol{\cal {P}}_{i}
}\longrightarrow (R \times R)
     _{\overline P_i}\,,\;\;\; \mbox{for $1 \leq i \leq 2$.}
  $$
\item[(b)] The restriction of $\gamma$
$$
\gamma \left|_{V_Y(\boldsymbol{\mathfrak O_1})\cap
V_Y(\boldsymbol{\mathfrak O_2})}:
   V_Y(\boldsymbol{\mathfrak O_1})\cap
V_Y(\boldsymbol{\mathfrak O_2}) \ \longrightarrow \
V_X(I)\right.
$$
is a scheme isomorphism.

\item[(c)]  If $P\in \Spec(R)$ is such that $P \supseteq I$
and $ {\boldsymbol{\cal P}}\in \Spec (\RI)$ is the unique prime
ideal such that ${\boldsymbol{\cal P}}\cap R=P$, the following
diagram of canonical homomorphisms:
$$
\begin{CD}
(\RI)_{  \boldsymbol{\cal P}}  @>>>    R_P \\
@VVV       @Vu_PVV  \\
R_P\times R_P   @>v_P>> R_P\times (R_P/I_P)
\end{CD}
$$
is a pullback (where $I_P  := IR_P$, $u_P(x) :=(x,x+I_P)$ and
$v_P((x,y)):= (x,y+I_P)$, for $x, y \in R_P$ ), i.e. $(\RI)_{
\boldsymbol{\cal P}} \cong R_P\!\Join\!I_P$ (Proposition
\ref{pullback-I}).
\end{description}
\end{thm}

\noindent{\bf Proof.}   (a) Since  $\boldsymbol{\mathfrak
O_1}=\{0\}\times I$ (respectively,  $\boldsymbol{\mathfrak
O_2}=I \times \{0\}$)   is a common ideal of $R\times R$ and
$\RI$, this statement follows   from the general results on
pullbacks \cite[Theorem 1.4]{F}  and   from
Theorem \ref{spec2} (and its proof).  
% and from the fact that the canonical
%ring homomorphisms $R_P \hookrightarrow (R\times R)_{\overline
%P_i} $ are isomorphisms, for $1 \leq i \leq 2$.
 Note that
$Z\setminus V_Z(\boldsymbol{\mathfrak O_1}) \cong
 \left( (X \amalg X)\setminus (X \amalg V_{X}(I))\right)=
  X \setminus V_{X}(I)=\left( (X \amalg X)\setminus (V_X(I)\amalg X)\right)
  \cong Z\setminus V_Z(\boldsymbol{\mathfrak O_2})$.

\smallskip
\noindent (b) Note that $V_Y(\boldsymbol{\mathfrak O_1}) \cap
V_Y(\boldsymbol{\mathfrak O_2}) = V_Y(\boldsymbol{\mathfrak
O_1}+\boldsymbol{\mathfrak O_2})$ and $\boldsymbol{\mathfrak
O_1}+\boldsymbol{\mathfrak O_2} = I \times I$. Therefore the
present statement follows from the fact that the canonical
surjective homomorphism $\RI \rightarrow R/I$, defined by $(r,
r+i) \mapsto r+I$ (for each $r \in R$ and $i \in I$) has kernel
equal to $I \times I$.
%The first part of the statement follows by Theorem
%\ref{spec2} (part (b)): $V_Y(\boldsymbol{\mathfrak O_1}) \cap
%V_Y(\boldsymbol{\mathfrak O_2})=V(I\times I)$ and every prime
%ideal of $\RI$ containing $I\times I$ is of the form $(P\times R)
%\cap (\RI)= (R\times P)\cap(\RI)=\boldsymbol(\mathcal P)$, for
%some $P\ in V_X(I)$. Moreover, by Proposition \ref{spec},
%$(\RI)/\boldsymbol(\mathcal P)\cong R/P$.

\noindent (c)  If we start from the pullback diagram
considered in Proposition \ref{pullback-I} and we apply the tensor
product $R_P \otimes_R \mbox{---}$\,, \ then by \cite[Proposition
1.9]{F} we get the following pullback diagram:
$$
\begin{CD}
R_P\otimes_R (\RI)   @>   id \otimes v' >>   R_P \otimes_R R  \\
@V  id \otimes u' VV       @V  id \otimes u  VV  \\
R_P \otimes_R (R\times R)   @>  id \otimes v  >> R_P \otimes_R
(R\times (R/I))\,.
\end{CD}
$$
  Note that, by the properties of the tensor product, we deduce
immediately the following canonical ring isomorphisms: $R_P \otimes_R
(R\times R) \cong R_P \times R_P$, $R_P \otimes_R R \cong R_P$ and
that $R_P \otimes_R (R\times (R/I))\cong R_P \times (R_P \otimes_R
(R/I)) \cong R_P \times (R_P/IR_P)$. Therefore, the previous pullback
diagram gives rise to the following pullback of canonical
homomorphisms:
$$
\begin{CD}
R_P\otimes_R (\RI)   @> >>   R_P\\
@V VV       @V u_P VV  \\
R_P \times R_P  @> v_P  >>  R_P \times (R_P/I_P)\,.
\end{CD}
$$
On the other hand, recall that $\Spec(R_P \otimes_R (\RI))$ can be
canonically identified (under the canonical homeomorphism
associated to the natural ring homomorphism $\RI \rightarrow R_P
\otimes_R (\RI) $) with the set of all prime ideals
${\boldsymbol{\cal H}} \in \Spec(\RI)$ such that
${\boldsymbol{\cal H}} \cap R \subseteq P$.  Since we know already
that,  in the    present situation,   there exists a unique prime
ideal $  {\boldsymbol{\cal P}}  \in \text{Spec}(\RI)$ such that $
{\boldsymbol{\cal P}}  \cap R=P$ (Theorem \ref{spec2}  (1, b))   and
that the canonical embedding $R \hookrightarrow \RI$ has  the
going-up property, we deduce that $\Spec(R_P \otimes_R (\RI))$ can
be canonically identified with the set of all the prime ideals of
$\RI$ contained in $  {\boldsymbol{\cal P}}$. Therefore $R_P
\otimes_R (\RI)$ is a local ring with a unique maximal ideal
corresponding to the prime ideal   $ {\boldsymbol{\cal P}}$ of
$\RI$ and thus we deduce that the canonical ring homomorphism
$(\RI)_{\boldsymbol{\cal P}} \rightarrow R_P \otimes_R (\RI)$ is
an isomorphism. \hfill $\Box$

%PROPOSITION 3.9
\begin{prop}\label{pullback-II}
The ring $\RI$ can be obtained as a pullback of the following
diagram of canonical homomorphisms:
$$
\begin{CD}
R \Join I  @> \widetilde{v}'>>    R/I \\
@V \widetilde{u}'VV       @V\widetilde{u}VV  \\
R\times R   @>\widetilde{v}>> R/I\times R/I
\end{CD}
$$
where $\widetilde{u}$ is the diagonal embedding, $\widetilde{v}$
is the canonical surjection $(x,y) \mapsto (x+I, y+I)$,
$\widetilde{u}'$  is the natural inclusion and $\widetilde{v}'$ is
defined by $(x, x+i) \mapsto x+I$, for all $x, y \in R$ and $i \in
I$.
\end{prop}

\noindent   {\bf Proof.}  By Proposition \ref{pullback-I} we know
that
$$
\begin{CD}
R \Join  I   @>>>   R \\
@VVV       @VuVV  \\
R\times R   @>v>> R\times R/I
\end{CD}
$$
is a pullback.  On the other hand, it is easy to verify that the
following diagram:
$$
\begin{CD}
  R  @>   \varphi >>    R/I \\
@V  u  VV       @V\widetilde{u}VV  \\
R\times R/I   @>w>> R/I\times R/I
\end{CD}
$$
is a pullback, where $w$ is the canonical surjection $(x,y)
\mapsto (x+I,y)$   and $\varphi$ is the natural proiection $x
\mapsto x +I$, for each $x \in R$ and for each $y \in R/I$.  The
conclusion follows by   juxtaposing  two pullbacks. \hfill $\Box$

%COROLLARY 3.10 {pullback-M}
\begin{cor} \label{pullback-M}
If $R$ is a local ring, integrally closed in $T(R)$ with maximal ideal $M$ and residue field
$k$,  then $\RM$ is
seminormal in its integral closure inside $T(R) \times T(R)$
(which, in this situation, coincides with  $R \times R$).
\end{cor}

\noindent  {\bf Proof.} By the previous proposition $R \Join M$
(which is a local ring) can be obtained as a pullback of the
following diagram of canonical homomorphisms:
$$
\begin{CD}
R \Join M  @>   \widetilde{v}' >>    k \\
@V  \widetilde{u}' VV       @V\widetilde{u}VV  \\
R\times R   @>\widetilde{v}>> k\times k
\end{CD}
$$
The statement follows from the fact that, in this case, the
integral closure of $\RM$ in  $T(R) \times T(R)$ coincides with $R
\times R$ (Corollary \ref{noeth+intclos} (c)).  Therefore, since
$\widetilde{u}$ is a minimal extension, then  $\widetilde{u}' $ is
also minimal \cite[Lemme 1.4  (ii)]{f-o}, and thus   the
conclusion follows from  \cite[Th\'eor\`eme 2.2 (ii))]{f-o} and
from \cite[(1.1)]{t} (keeping in mind Theorem \ref{spec2} (c)).
\hfill $\Box$

%EXAMPLE 3.11
\begin{ex}
\rm {\bf (a)} Let $R :=  k[[t]]$ (where $k$ is a field and $t$ an
indeterminate) and let $I  := t^nR $. Using Proposition
\ref{pullback-II}, if we denote by $h^{(i)}(t)$ the  $i$--th
derivative of a power series $h(t) \in k[[t]]$,  it is easy to see
that
$$\RI=\{(f(t),g(t)) \mid f(t), g(t) \in R\,,\;   f^{(i)}(0)=g^{(i)}(0)\ \forall \
i=0,\dots n-1 \}\ .
$$

\smallskip
\noindent {\bf (b)} Let $R := k[x,y]$ and $I := xR$.  In this case
$$\RI=\{(f(x,y),g(x,y)) \mid  f(x,y), g(x,y)\in R\,,\; f(0,y)=g(0,y) \}\ .$$
Setting $Y=\Spec(\RI)$ and $X=\Spec(R)$, by Proposition
\ref{spec}, $V_Y(\boldsymbol{\mathfrak O_i}) \cong
 \Spec(k[x, y])$.   On the other hand, by Theorem \ref{loc},
$V_Y(\boldsymbol{\mathfrak O_1})\cap V_Y(\boldsymbol{\mathfrak
O_2})=  V_Y((xR\times xR))\cong V_X(xR) \cong
 \Spec(k[y])$.  Hence the ring $\RI$ is the coordinate ring of two
affine planes with a common line.  Note that we can present $\RI$ as
quotient of a polynomial ring in the following way: consider the
homomorphism $\lambda: k[x,y,z]\ \longrightarrow \ R\times R$, defined by
$\lambda(x):= (x,x)$, $ \lambda(y):=(y,y)$ and $ \lambda(z):= (0,x)$. It is not difficult
to see that Im$(\lambda)=\RI$ and $\Ker(\lambda)=(zx-z^2)k[x,y,z]$.
\end{ex}
\smallskip

\noindent \bf Acknowledgment \rm

\noindent The authors want to thank the referee for his/her many valuable suggestions and comments which have improved the final version of this paper.

\bigskip

%
%\noindent
%Authors' adresses:\\

\noindent
Marco D'Anna\\ Dipartimento di Matematica e Informatica\\
Universit\`a di Catania\\ Viale Andrea Doria 6\\ 95125 Catania,
Italy\\
{\texttt {mdanna@dipmat.unict.it}} \\

\smallskip\noindent
Marco Fontana\\ Dipartimento di Matematica\\
Universit\`a di Roma degli Studi ``Roma Tre''\\ Largo San Leonardo
Murialdo 1\\ 00146 Roma, Italy\\
\texttt{fontana@mat.uniroma3.it}


\begin{thebibliography}{99}

\footnotesize


\bibitem{d'a} M.~D'Anna,  {\em  A construction of Gorenstein
rings}, \rm J. Algebra (to appear).

 \bibitem{d'a-f} M.~D'Anna and M. Fontana, {\em The amalgamated duplication of a ring along a multiplicative-canonical ideal}, Preprint 2006.

\bibitem{f-o} D. Ferrand and J.-P. Olivier, \it Homomorphismes
mimimaux d'anneaux, \rm  J. Algebra \bf16 \rm (1970), 461--471.

\bibitem{F} M.~Fontana, {\em Topologically defined classes of
commutative rings},
Ann. Mat. Pura Appl.  \bf 123 \rm (1980), 331--355.

\bibitem{fs} R. Fossum,
\it Commutative extensions by canonical
modules are Gorenstein rings,
\rm Proc. Am. Math. Soc. \bf 40 \rm (1973), 395--400.

\bibitem{fgr} R. Fossum, P. Griffith and I. Reiten, \it  Trivial
extensions of Abelian categories. Homo\-lo\-gical algebra of trivial
extensions of Abelian categories with applications to ring theory,
\rm
Lecture Notes in Mathematics \bf 456\rm , Springer-Verlag, Berlin,
1975.

\bibitem{H-G} S. Gabelli and E.G. Houston, \it  Ideal theory in pullbacks, \rm in ``Non-Noetherian Commutative Ring Theory'', S.T. Chapman and S. Glaz Eds., Kluwer Academic Publishers, 2000, 199--227. 

\bibitem{Gilmer} R. Gilmer, Multiplicative ideal theory, M. Dekker,
New York, 1972.

\bibitem{gl} S. Glaz, Commutative  coherent  rings, Lecture Notes in Mathematics \bf 1321\rm,  Springer-Verlag, Berlin, 1989.

 \bibitem{hhp} W. Heinzer, J. Huckaba and I. Papick, \it m--canonical
ideals in integral domains, \rm Comm. Algebra \bf 26 \rm (1998),
3021--3043.

\bibitem{h} J. Huckaba, Commutative rings with zero divisors, M.
Dekker, New York, 1988.

\bibitem{Ka} I. Kaplansky, Commutative Rings, Allyn and Bacon,
Boston, 1970.

\bibitem {La} J. Lambek, Lectures on Rings and Modules, Blaisdell
Publishing Company, Waltham, 1966.


\bibitem{Matsumura} H. Matsumura,  Commutative ring theory, Cambridge
University Press, Cambridge, 1986.

\bibitem{n-1955} M. Nagata,
\it The theory of multiplicity in general local rings,
\rm Proc. Intern. Symp. Tokyo-Nikko 1955, Sci. Council of Japan,
Tokyo 1956, 191--226.

\bibitem{n}
M. Nagata,  {Local Rings}, Interscience, New York, 1962.


\bibitem{r} I. Reiten, \it The converse of a theorem of Sharp on
Gorenstein modules, \rm  Proc. Amer. Math. Soc. \bf 32 \rm (1972),
417--420.

\bibitem{t} C. Traverso,  \it  Seminormality and Picard group, \rm
Ann. Sc. Norm. Sup. Pisa \bf 24 \rm (1970), 585--595.


\end{thebibliography}
\end{document}